\documentclass[12pt]{article}
\usepackage{amsmath}
\usepackage{amssymb}
\numberwithin{equation}{section}

\title{Expectations of hook products on large partitions
}

\author{
M. Adler\thanks{2000 {\em Mathematics Subject
classification.} 60C55, 60F05, 05A05, 05A15, 34M55
.%
~\newline
{\em Key words and phrases.}  Integrals over
Grassmannians, longest increasing sequences in random
words, limit theorems, chi-square distribution,
Painlev\'e equations.
\newline  Department of Mathematics, Brandeis University,
Waltham, Mass 02454, USA. E-mail: adler@brandeis.edu.
The support of a National Science Foundation grant \#
DMS-01-00782 is gratefully acknowledged.}~~~~~~A.
~Borodin\thanks{Department of Mathematics, California
Institute of technology, Pasadena, CA 91125, USA and
Clay mathematics institute Long-term Prize Fellow.
E-mail: borodin@caltech.edu. }~~~~~ P. van
Moerbeke\thanks{ Department of Mathematics, Universit\'e
de Louvain, 1348 Louvain-la-Neuve, Belgium and Brandeis
University, Waltham, Mass 02454, USA. E-mail:
vanmoerbeke@math.ucl.ac.be and @brandeis.edu. The
support of a National Science Foundation grant \#
DMS-01-00782, a Nato, a FNRS and a Francqui Foundation
grant is gratefully acknowledged. This work was done
while PvM was a member of the Clay Mathematics
Institute, One Bow Street, Cambridge, MA 02138, USA.
 MA
and PvM would like to thank Persi Diaconis for
interesting conversations regarding section 4.}}

 \date{}

\catcode`\@=11
\let\c@equation=\relax
\newcounter{equation}[subsection]

\catcode`\@=12

\newcommand{\MAT}[1]{\left(\begin{array}{*#1c}}
\newcommand{\mat}{\end{array}\right)}
\newcommand{\qed}{\leavevmode\unskip\nobreak\penalty200\hskip2pt\null
\nobreak\hfill\rule{1.1ex}{1.1ex}%\parfillskip=0pt
\medbreak }

\newcommand{\rg}{\rightarrow}

\newcommand{\HR}{{\cal H}}

\newcommand{\LR}{{\cal L}}

\newcommand{\BC}{{\mathbb C}}

\newcommand{\BY}{{\mathbb Y}}

\newcommand{\iy}{\infty}
\newcommand{\pl}{\partial}
\newcommand{\al}{\alpha}

\newcommand{\vr}{\varepsilon}

\newenvironment{proof}
{\medskip\noindent{\it Proof:\/} }{\qed}
\newenvironment
        {remark}{\medskip\noindent\underline{\it Remark:\/} }{\medbreak}

\newcommand{\vp}{\varphi}
\newcommand{\la}{\langle}
\newcommand{\ra}{\rangle}

\newcommand{\dt}{\delta}
\newcommand{\Dt}{\Delta}

\newcommand{\BR}{{\mathbb R}}
\newcommand{\lb}{\lambda}

\newcommand{\no}{\nonumber}

\newcommand{\BJ}{{\mathbb J}}

\def\be#1\ee{\begin{equation}#1\end{equation}}
\def\bea#1\eea{\begin{eqnarray}#1\end{eqnarray}}
\def\bean#1\eean{\begin{eqnarray*}#1\end{eqnarray*}}

\newcommand{\Tr}{\operatorname{\rm Tr}}

\newtheorem{definition}{Definition}[section]
\newtheorem{theorem}[definition]{Theorem}
\newtheorem{lemma}[definition]{Lemma}
\newtheorem{corollary}[definition]{Corollary}
\newtheorem{proposition}[definition]{Proposition}

\makeatletter
\def\ps@X{\let\@mkboth\@gobbletwo
        \def\@oddhead{\tt %Adler-van\,Moerbeke:%
        \hfil\S\thesection,
        p.\thepage}
        \def\@oddfoot{\rm\hfil\thepage\hfil}
        \let\@evenhead\@oddhead
        \let\@evenfoot\@oddfoot}
\makeatother

\catcode `!=11

\newdimen\squaresize
\newdimen\thickness
\newdimen\Thickness
\newdimen\ll! \newdimen \uu! \newdimen\dd! \newdimen \rr! \newdimen
\temp!

\def\sq!#1#2#3#4#5{%
\ll!=#1 \uu!=#2 \dd!=#3 \rr!=#4
\setbox0=\hbox{%
 \temp!=\squaresize\advance\temp! by .5\uu!
 \rlap{\kern -.5\ll!
 \vbox{\hrule height \temp! width#1 depth .5\dd!}}%
 \temp!=\squaresize\advance\temp! by -.5\uu!
 \rlap{\raise\temp!
 \vbox{\hrule height #2 width \squaresize}}%
 \rlap{\raise -.5\dd!
 \vbox{\hrule height #3 width \squaresize}}%
 \temp!=\squaresize\advance\temp! by .5\uu!
 \rlap{\kern \squaresize \kern-.5\rr!
 \vbox{\hrule height \temp! width#4 depth .5\dd!}}%
 \rlap{\kern .5\squaresize\raise .5\squaresize
 \vbox to 0pt{\vss\hbox to 0pt{\hss $#5$\hss}\vss}}%
}%end of \hbox
 \ht0=0pt \dp0=0pt \box0
}%end of \sq!

\def\vsq!#1#2#3#4#5\endvsq!{\vbox to \squaresize{\hrule
width\squaresize height 0pt%
\vss\sq!{#1}{#2}{#3}{#4}{#5}}}

\newdimen \LL! \newdimen \UU! \newdimen \DD! \newdimen \RR!

\def\vvsq!{\futurelet\next\vvvsq!}
\def\vvvsq!{\relax
  \ifx     \next l\LL!=\Thickness \let\continue=\skipnexttoken!
  \else\ifx\next u\UU!=\Thickness \let\continue=\skipnexttoken!
  \else\ifx\next d\DD!=\Thickness \let\continue=\skipnexttoken!
  \else\ifx\next r\RR!=\Thickness \let\continue=\skipnexttoken!
  \else\def\continue{\vsq!\LL!\UU!\DD!\RR!}%
  \fi\fi\fi\fi
  \continue}

\def\skipnexttoken!#1{\vvsq!}

\def\place#1#2#3{\vbox to 0pt{\vss
\rlap{\kern#1\squaresize
  \raise#2\squaresize\hbox{$#3$}}
\vss}}

\def\Young#1{\LL!=\thickness \UU!=\thickness \DD! = \thickness \RR! =
\thickness \vbox{\smallskip\offinterlineskip \halign{&\vvsq! ##
\endvsq!\cr #1}}}

\catcode `!=12

\begin{document}

%vspace{-2cm}

 \maketitle

\begin{abstract}

Given uniform probability on words of length $M=Np+k$,
from an alphabet of size $p$, consider the probability
that a word ({\bf i}) contains a subsequence of letters
$p, p-1,\ldots,1$ in that order and ({\bf ii}) that the
maximal length of the disjoint union of $p-1$ increasing
subsequences of the word is $\leq M-N$. A generating
function for this probability has the form of an
integral over the Grassmannian of $p$-planes in $\BC^n$.
The present paper shows that the asymptotics of this
probability, when $N\rightarrow \iy$, is related to the
$k^{\mbox{\tiny th}}$ moment of the
$\chi^2$-distribution of parameter $2p^2$. This is
related to the behavior of the integral over the
Grassmannian $Gr(p,\BC^n)$ of p-planes in $\BC^n$, when
the dimension of the ambient space $\BC^n$ becomes very
large. A different scaling limit for the Poissonized
probability is related to a new matrix integral, itself
a solution of the Painlev\'e IV equation. This is part
of a more general set-up related to the Painlev\'e V
equation.

\end{abstract}

% \newpage

\tableofcontents

\section{Introduction}

%\vspace{1cm}
Consider the set of words
  $$
  \pi\in S_{\ell}^p:=
\left\{\mbox{words $\pi$ of length $\ell$, built from an
alphabet $\{1,...,p\}$}\right\},
 $$
  (with $|S_{\ell}^p|=p^{\ell}$) with the uniform probability distribution
\be
 P^{\ell,p}(\pi)=\frac{1}{p^{\ell}}.
  \ee
Let $\BY$ denote the set of all partitions $\lb$. Let
 $\lb^{\top}$ be the dual partition, i.e.,
  obtained by flipping the Young diagram $\lb$ about
  its diagonal. So, $\lb_1^{\top}$ is the length of the
  first column of $\lb$.
  Moreover $h^{\lb}$ denotes the product of the
  hook lengths
    $
  h^{\lb}_{ij}:=\lb_i+\lb_j^{\top}-i-j+1
   $
 and also $1^p$ denotes the infinite vector
 $$
1^p:=(\overbrace{1,\ldots,1}^p,0,0,\ldots).$$

 The RSK correspondence between words and
 pairs of semi-standard and standard tableaux
 induces a probability measure on partitions
  \be
  \lb \in \BY_{\ell}=\{\mbox{partitions
 $ \lb \in \BY$ of weight $|\lb|=\ell$} \},
  \ee
 given by
    \be
  P^{\ell,p}(\lb)=
        \frac{f^{\lb}~{\bf s}_{\lb}(1^p)}{p^{|\lb|}}
 ,\label{probability on words}
  \ee
 %
%where $h:=(h_1,...,h_p)$ with $h_i:=p+\lb_i-i$, upon
%using (\ref{standard}) and (\ref{semi-standard}).
 having
$$(\mbox{ support }P^{\ell,p})\subseteq
 \BY_{\ell}^{(p)}:=\{\lb \in \BY_{\ell},~ \mbox{such that}~
\lb_1^{\top}\leq p \}.$$
The symbol ${\bf s}_{\lb}$ denotes the Schur polynomial
associated
 with the partition $\lb$.
 Besides the probability $P^{\ell,p}$, we also consider the corresponding
Poissonized measure, depending on the real parameter
$x$,
\be
 P_{x,p}(\lb)=e^{-px}\frac{(xp)^{\vert\lb\vert}}
  {\vert\lb\vert
!}P^{\vert\lb\vert,p}(\lb),\qquad\lb\in\BY^{(p)},
\label{Poissonized probability}\ee
 having
$$(\mbox{ support }P_{x,p})\subseteq
 \BY_{ }^{(p)}:=\{\lb \in \BY_{ },~ \mbox{such that}~
\lb_1^{\top}\leq p \}.$$

We discuss very briefly the combinatorics needed in this
problem; for more details, see
\cite{MacDonald,Sagan,Stanley}.
For the partition $\lb$, define the symbol
 \be
 (n)_{\lb}:=\prod_{i}(n+1-i)_{\lb_i},\mbox{ with }
 (x)_n=x(x+1)\ldots (x+n-1),~x_0=1
  %n(n+1)\ldots(n+\lb_i-1),
  \label{symbol}
  \ee
Define for $q\geq \lb_1^{\top}$, (throughout $ \Dt_q$
denotes the Vandermonde determinant in $q$ variables)
 \bea
   f^{\lb}=\# \left\{\begin{array}{l}
 \mbox{standard tableaux}\\ \mbox{of shape
$\lb$}\end{array}\right\}
 &=&\frac{|\lb
|!}{h^{\lb}}=
  \frac{|\lb|!}{u^{|\lb|}}~\left.  {\bf s}_{\lb}(x)
   \right|_{\sum_{\ell}x_{\ell}^i=\delta_{1i}u}
   \nonumber\\
  &=&  |\lb|! ~
  \frac{\Delta_{q}(q+\lb_1-1,\ldots,q
 +\lb_{q}-q)}
{\displaystyle{\prod_1^{q}} (q+\lb_i-i)!},
 \nonumber \\
  \label{standard}\eea
and
\bea \# \left\{\begin{array}{l}\mbox{semi-standard
tableaux
}\\
  \mbox{of shape $\lb$ filled with} \\ \mbox{numbers
from $1$ to $q$}
\end{array}
\right\}
 =
{\bf s}_{\lb}(1^q)
  \!&=&\!
    \displaystyle{
\frac{\Delta_q(q\!+\!\!\lb_1\!-\!1,\ldots,q\!+\!\lb_q\!-\
\!q)}
{\displaystyle{\prod^{q-1}_{i=1}}i!}} \no\\
 \!&=\!&\!
 \prod_{(i,j)\in \lb}\frac{j-i+q}{h^{\lb}_{ij}}
 \no\\& =&  \frac{(q)_{\lb}}{h^{\lb}}
   \label{semi-standard}\eea

 A subsequence $\sigma$ of the word
$\pi$ is {\em weakly
 $k$-increasing}, if it can be written as
 \be
 \sigma=\sigma_1\cup \sigma_2\cup \ldots \cup
 \sigma_k,
  \ee
 where $\sigma_i$ are disjoint weakly increasing subsequences of the word
 $\pi$, i.e., possibly with repetitions.
  The length of the longest increasing/decreasing
  subsequences is closely related to the shape of the
  associated partition, via the RSK correspondence
  :
 \bea
 d_1(\pi)&=& \left\{\begin{array}{l} \mbox{length of
the longest {\em strictly}}\\\mbox{ decreasing
 subsequence of $\pi$}
\end{array}
\right\}=\lambda^{\top}_1
 \nonumber\\
i_k(\pi)&=&
 \left\{\begin{array}{l}
\mbox{length of the
 longest {\em weakly}}\\ \mbox{ $k$-increasing
 subsequence of $\pi$}
\end{array}
\right\}= \lambda_1+\ldots+\lambda_k
 \nonumber\\
 \label{RSK}\eea
Define the generalized hypergeometric function in terms
 of the symbol (\ref{symbol}), viewed as a symmetric function
 in an infinite number of variables $x_i$
\be {}_2F_1^{(1)}(p,q;n;x): =\sum_{\kappa \in \BY}
 \frac{(p)^{ }_{\kappa}(q)^{ }_{\kappa}}
  {(n)^{ }_{\kappa}}
~  \frac{{\bf s}_{\kappa}(x)}
 {h^{\kappa}},
 \ee
 which, upon restriction, using power sums and upon using (\ref{standard}),
 yields
%  $$
%  \left.  {\bf s}_{\lb}(x)
%   \right|_{\sum_{\ell}x_{\ell}^i=\delta_{1i}u}
%   = \frac{u^{|\lb|}}{h^{\lb}}
%   ,$$ one finds, upon restriction,
\bea
   {}_2F_1^{(1)}(p,q;n;x)
 \Bigr|_{\sum_{\ell}x_{\ell}^i= \delta_{1i}u}
&=&
 \sum_{\kappa \in \BY}u^{|\kappa|}
 \frac{(p)^{ }_{\kappa}(q)^{ }_{\kappa}}{({h^{\kappa}})^2
  (n)^{ }_{\kappa}}
 \no \\
 &=&
  \sum_{k=0}^{\iy}u^k \sum_{\kappa \in \BY_k}%u^{|\kappa|}
 \frac{(p)^{ }_{\kappa}(q)^{ }_{\kappa}}{({h^{\kappa}})^2
  (n)^{ }_{\kappa}}
 \label{restiction of hypergeometric}\eea

 As a reminder, the chi-square distribution of
parameter $m$ is the distribution of
$$
Z_m=\sum_1^mX_i^2 ,$$
 where the $X_i$'s
 are $m$ independent normal $N(0,1)$-random variables.
%Theorem 1.1 and Corollary 1.2 involve the hook length,
%

\vspace{1cm}
 %  \newpage

   $\overbrace{\hspace{9.4cm}}^{\mbox{\tiny
partition}~\lb}$

$\overbrace{\hspace{6cm}}^{\mbox{\tiny partition}~\mu}$

\bigbreak

\hspace{-.8cm}$ p\left\{ \begin{array}{c} \\ \\ \\ \\
\end{array}\right.
 $

%$\overbrace{10cm}$
 \vspace{-2.2cm}
$ \squaresize .3cm \thickness .01cm \Thickness .07cm
\Young{
 &&&&&&&&&&&&&&&r&&&&r&&&&&&&&&&&\cr
 &&&&&&&&&&&&&&&r&&&&r&&&&&&&&&\cr
 &&&&&&&&&&&&&&&r&&&&r&&&&&&&\cr
 &&&&&&&&&&&&&&&r&&&&r&&&&&&\cr
 &&&&&&&&&&&&&&&r&&&&r&&&&&&             \cr
 &&&&&&&&&&&&&&&r&&&&r&&                 \cr
 &&&&&&&&&&&&&&&r&&&&r                   \cr }\hspace{.9cm}
   $

\vspace{-.8cm}

 $$%\hspace{.1cm}
  \underbrace{\hspace{4.5cm}}_{n-q}
%\hspace{.01cm} ~~~
~\underbrace{\hspace{.8cm}}_{q-p}
 \underbrace{\hspace{3.7cm}}_{\#\mbox{\tiny boxes}~=k}
  \hspace{2.5cm}$$

$$\mbox{Figure 1}$$

 The expectations $E_{x,p}$ and $E^{\ell,p}$ are taken
  with
 regard to the probabilities $P_{x,p}$ and $P^{\ell,p}$
 defined above. The functions on partitions, of which
 the expectations are taken, are products of hook
 lengths restricted to a vertical strip in the partition
 of width $q-p$, as in Figure 1.
In \cite{AvM3}, expectations of this type
%(\ref{Theorem 5.1})
have been studied and linked to integrals over the
Grassmannian space Gr(p,${\BC^n}$) of $p$-planes in
$\BC^n$, with regard to the Weyl measure $d\rho(Z)$;
besides, these integrals relate to specific solutions of
the Painlev\'e V equation.

We now state the following proposition, established in
\cite{AvM3}:
\begin{proposition}
 Given the partition
       \be\mu=(n-p)^p:=(\overbrace{n-p,n-p,\ldots,n-p}^p).\ee
%%For fixed $p\leq q\leq n/2$, the generating function for
the mathematical expectation of the product of the hook
lengths over the strip of width $q-p$, with regard to
the probability (\ref{probability on words}), is given
by\footnote{Remember $(\mbox{ support }P^{\ell,p})
 \subseteq\{\lb \in \BY,~ \mbox{such that}~
|\lb|=\ell,~\lb_1^{\top}\leq p \}.$}

\bea
 \lefteqn{e^{px }~E_{x,p}\left( I_{\{\lambda \supseteq \mu \}}(\lb)
  \prod_{{(i,j)\in \lb}
 \atop {n-q<j\leq n-p}}
  h_{(i,j)}^{\lb}
   \right)}\no\\
 &=&%\hspace{-1cm}
 \displaystyle{
%\frac{}
 \sum_{\ell\geq p(n-p)}
  \frac{(px)^{\ell}}{\ell!}
 E^{\ell,p}\left( I_{\{\lambda \supseteq \mu \}}(\lb)
  \prod_{{(i,j)\in \lb}
 \atop {n-q<j\leq n-p}}
  h_{(i,j)}^{\lb}
   \right) }
   \no
 %  \eea\bea
 \no \\&&\no\\
   &=&
   \tilde{\tilde c_{ }}
    x^{(n-p)p}~{}_2F_1^{(1)}(p,q;n;y)
 \Bigr|_{\sum_{\ell}y_{\ell}^i= \delta_{1i}x}
  \nonumber \\&=&
  {\tilde c_{ }} x^{(n-p)p}
\int_{Gr(p,\BC^n)}e^{x\Tr(I+Z^{\dag}Z)^{-1}}\det(Z^{\dag
}Z)^{-(q-p)} d\rho(Z)
 % \eea\bea\no%
 \no\\
 &=&
    \tilde{\tilde c_{ }}
    x^{(n-p)p}
    \exp {\displaystyle{\int_0^x\frac{u(y)-p(n-p)+py}{y}
    dy}
    },\label{proposition4.1}
 \eea
       where\footnote{with
\bean
 \tilde c_{ }^{-1} :=  \prod_{i=1}^p i!~(n-q-i)!
 ~~~~\mbox{and}~~~~
 \tilde{\tilde c_{ }} := \frac{(q-i)!}{(n-i)!}.\eean}
$u(x)$
  is the unique solution to the initial value problem:
 \be
 \left\{\begin{array}{l} \displaystyle{ x^2
u^{\prime\prime\prime}+xu^{\prime\prime}
 +6x{u^{\prime} }^2-4uu^{\prime}+4Qu^{\prime}
 -2Q^{\prime}u+2R=0 } \\ \mbox{with} \hspace{9cm}
  \mbox{\bf (Painlev\'e V)} \\
\displaystyle{ u_{}(x)=p(n-p)-
 \frac{p(n\!-\!q)}{n}x+\ldots+
   a_{n+1}x^{n+1}+O(x^{n+2})  ,~\mbox{near}~ x=0
,}
\end{array}
 \right.
 \label{1.0.14}\ee
  with a specific coefficient $a_{n+1}$ and where
  \bea
 4Q&=&
 -x^2+2(n+2(p-q))x-(n-2p)^2 \nonumber\\
 2R&=&
  p(p-q)(x+n-2p).\label{QR-polynomials}
    \eea

\end{proposition}

\vspace{1cm}

 This paper is concerned with what happens
when $n \rightarrow \iy$. This is related to the
behavior of the integral over the Grassmannian
$Gr(p,\BC^n)$ of p-planes in $\BC^n$, when the dimension
of the ambient space $\BC^n$ becomes very large.
To be precise:%We now state:

\begin{theorem} Given the partition $\mu=(n-p)^p$, as in Figure 1, the following expectation behaves, for
large $n$, like the moments of the chi-square
distribution of parameter $2pq$ :
\bea
  \lefteqn{
\lim_{n\rg\iy}n^{\frac{p^2-1}{2}}E^{p(n-p)+k,p}
  \left(
I_{\lb\supseteq\mu}(\lb)\prod_{n-q<j\leq
n-p}h^{\lb}_{ij}\right)
 }\no\\
 %  \no\\
  &=&
  \frac{
   \sqrt{p}
 \displaystyle{ \prod^p_1}(q-j)!}
   {(\sqrt{2\pi})^{p-1}k!}
E\left(\left( \frac{1}{2}Z_{2pq}\right)^k\right),
 \label{Theorem 5.1}\eea
and so the expectation decays, when $n\rg \iy$, as
$$
E^{k+p(n-p),p}
  \left(
I_{\lb\supseteq\mu}(\lb)\prod_{n-q<j\leq
n-p}h^{\lb}_{ij}\right) \simeq
  c_{p,q,k}~
  n^{-\frac{p^2-1}{2}}
%   \frac{
%   \sqrt{p}
% \displaystyle{ \prod^p_1}(q-j)! }
%   {(\sqrt{2\pi})^{p-1}  }
% \left( {{pq-1+k}\atop k}\right)
%
%
,$$
 with
 $$
 c_{p,q,k}:=  \frac{
   \sqrt{p}
 \displaystyle{ \prod^p_1}(q-j)! }
   {(\sqrt{2\pi})^{p-1}  }
 \left( {{pq-1+k}\atop k}\right)
 .$$

\end{theorem}

The next statement deals with the special case, where
$q=p$. Namely, setting $N:= n-p$, what is,
asymptotically for large N, the probability
$P^{Np+k,p}(\lambda_p\ge N)$ ?
\begin{corollary} Given an alphabet of size $p$
and an integer $k>0$, we give the behavior of the
probability on the set of words of length $Np+k$ for
large $N$. Notice $i_p(\pi)=\{$length of the
word$\}=Np+k$ automatically, when $d_1(\pi)=p$. So,
$i_{p-1}(\pi)$ is the first non-trivial quantity. We now
have:
   \bean
 \lefteqn{\hspace{-1.8cm}
 \lim_{N\rg\iy}
 N^{\frac{p^2-1}{2}}
 P^{Np+k,p}(\lambda_p\ge N)}\\
 &=&\lim_{N\rg\iy}
 N^{\frac{p^2-1}{2}}
 P^{Np+k,p}\left(
  \begin{array}{c}
     d_1(\pi)=p\\  \\
     i_{p-1}(\pi)\leq N(p-1)+k
     \end{array}
     \right) \\
     & =&
  \frac{
   \sqrt{p}
 \displaystyle{ \prod^p_1}(p-j)!}
   {(\sqrt{2\pi})^{p-1}k!}
   E\left(\left( \frac{1}{2}Z_{2p^2}\right)^k\right).
   \eean
 Thus the following decay holds
  for $N\nearrow \iy$,
 $$
 P^{Np+k,p}\left(
  \begin{array}{c}
     d_1(\pi)=p\\  \\
     i_{p-1}(\pi)\leq N(p-1)+k \end{array}
     \right)
 \simeq
 c_{p,p,k}
 \left(\frac{1}{N}\right)^{\frac{p^2-1}{2}}
%   \frac{ \sqrt{p}
% \displaystyle{ \prod^p_1}(p-j)!  }
%   {(\sqrt{2\pi})^{p-1}  }
%    \left( {{p^2-1+k}\atop k}
%    \right)
%
.$$

\end{corollary}

\bigbreak

The proofs of Theorem 1.2 and Corollary 1.3 will be
given in section 4.
In the next statement, we consider the expectation for
the Poissonized probability $P_{x,p}$ on Young diagrams
$\lb\in\BY_{\ell}^{(p)}$, as defined in
(\ref{Poissonized probability}). Define
 $$ \HR_p=\{p\times p\mbox{~Hermitian
matrices}\},
$$
and, for an interval $I\subset \BR$,
$$ \HR_p(I)=\{M\in \HR_p \mbox{~with spectrum in $I$}\}.
$$

\bigbreak

\begin{theorem}
Take $x>0,~ s\in \BR$ and set
% by means of %etting $n$ grow like %For $\mu=(n-p)^p $,
% and
%
\be n-p=x+s\sqrt{2x}, %~\mbox{for $x$ large,}
\label{rescaling} \ee
 we have, upon expressing $n$ in terms of $x$ by means of the rescaling
 (\ref{rescaling}) \footnote{with $\hat c=\tilde c
 2^{\frac{p(p-q)}{2}}$. Remember $\tilde c$ from footnote 2.}
\bea
\lefteqn{\lim_{x\rg\iy}\frac{1}{(2x)^{\frac{p(q-p)}{2}}}
E_{x,p}\left(I_{\lb\supseteq\mu}(\lb)\prod_{(i,j)\in
\lb\atop{n-q<j\leq n-p}}
h^{\lb}_{ij}\right)}\nonumber\\
& & \nonumber\\
 &=&
  \lim_{x\rg\iy} {\hat c} e^{-px}
    x^{p(n-\frac{p+q}{2})}
     \int_{Gr(p,\BC^n)}
     e^{x\Tr(I+Z^{\dag}Z)^{-1}}\det(Z^{\dag
}Z)^{-(q-p)} d\rho(Z)
 \no\\&=& p!\,\frac
 {\displaystyle{\int_{\HR_p [s, \iy)} \det (M-sI)^{q-p}e^{-\Tr M^2}dM}}
 {\displaystyle{\int_{\HR_p}  e^{-\Tr M^2}dM}}
 \nonumber\\
 &=&
 c ~\exp{\displaystyle \int_0^s h(y)dy} \label{main theorem}
 ,\eea
 %&=&\mbox{constant~}E_{\HR_p}\prod^p_{i=1}\left((z_i-c)^{q-p}I_{[c,\iy)}
%(z_i)\right)
%
 with $h(y)$  satisfying the Painlev\'e IV
  equation:
 \be
h^{\prime\prime\prime}+6h^{\prime
2}-4(y^2+2(q-2p))h^{\prime}+4yh-8(q-p)p =0.\label{ode in
h-1}
 \ee

\end{theorem}

\noindent The proof of Theorem 1.3 will be given in
sections 3 and 6. In section 2 we show that the
logarithmic derivative of a general multiple integral,
involving the square of a Vandermonde, satisfies a third
order differential equation, from which we derive, in
section 3, that the logarithmic derivative of the matrix
integral in (\ref{main theorem}) satisfies the
Painlev\'e equation. Section 6 contains a discussion on
how the Painlev\'e IV equation can be obtained from
Painlev\'e V, by means of the rescaling
(\ref{rescaling}), for large $x$. The following result
shows that a certain general integral satisfies
Painlev\'e V:

\begin{theorem}
Given the weight
$$
\rho(z)=(z-a)^{\al} (b-z)^{\beta}e^{\gamma z} ,\mbox{
with  } \al, \beta >-1 ,$$
on the interval $[a,b]$,
 the integral
  \be g(x):=\frac{\pl}{\pl
x}\log\int_{[a,b]^n}\Dt^2(z) \prod_1^n
e^{xz_k}\rho(z_k)dz_k
%\label{integral g}
\ee
 is a solution to the Painlev\'e V equation. To be
 precise,
\be
f(y)=n(n+\al+\beta)-\frac{y}{b-a}
    \left(g\Bigl(\frac{y}{b-a}-\gamma\Bigr)
   -na\right)
\label{f(y)}\ee
satisfies a version of the Painlev\'e V equation,
\bea &f^{ \prime\prime 2}& +\frac{4}{P^2}
 \left( (Pf^{\prime 2}+Q f^{\prime}+R)f^{\prime}
 - (P' f^{\prime 2}+\frac{}{}Q' f^{\prime}+R')f^{}
  \right. \nonumber\\ && \hspace{1.3cm} \left.
   +\frac{1}{2}(P^{\prime\prime}f^{\prime
}+Q^{\prime\prime} )f^2
 -\frac{1}{6} P^{\prime\prime\prime}f^3
   -\frac{1}{4} \beta^2n^2  \right)=0;
\label{1.0.21}\eea with %$\al$ an integration constant and
\bean
 P(y)&=&y\\
    4Q(y)&=&-y^2+2y(2n+\al-\beta)-(\al+\beta)^2 \\
    2R(y)&=& -\beta n(\al+\beta+y)
 .\eean

\end{theorem}

Note that (\ref{1.0.21}) is a well known form of
Painlev\'e V, as discussed in the Appendix. Theorem 1.5
and also Theorem 3.1 below (which is exactly the last
equality in (\ref{main theorem}), can be derived form
the work of Forrester and Witte \cite{Forrester
Witte1,Forrester Witte2}.

%\newpage

\section{A differential equation for a matrix integral }

\begin{proposition} Consider a weight $\rho(z)$
on an interval $E\subseteq\BR$, with rational
logarithmic derivative of the form
 \be
-\frac{\rho'(z)}{\rho(z)}=\frac{b_0+b_1z+b_2z^2}{a_0+a_1z+a_2z^2}
=: \frac{B(z)}{A(z)},
 \label{rational logderivative}
 \ee
and boundary condition
 \be
  A(z)\rho(z)z^k\bigl|_{\pl
E}=0,~~\mbox{for all $k=0,1,2,\ldots$.}
 \label{boundary condition}\ee
The expression\footnote{When $E$ is a finite interval,
the integral always converges, and for a infinite
interval, one may have to require $x>\alpha$ or $x<
\alpha$, for some $\alpha\in \BR$.
 }
 \be g(x):=\frac{\pl}{\pl
x}\log\int_{E^n}\Dt^2(z) \prod_1^n e^{xz_k}\rho(z_k)dz_k
\label{integral g} \ee
satisfies a third order differential equation,
equivalent to Painlev\'e V, namely:
 \be
 g^{\prime\prime\prime}+6g^{\prime
2}+\frac{4a_2(g^{\prime\prime}
+2gg^{\prime})}{a_2x-b_2}
 +\frac{2a_2^2g^2+P_2g'}{(a_2x-b_2)^2}
  +\frac{P_1 g -nQ_1}{(a_2x-b_2)^3}=0
 ,\label{third order ode}
 \ee
which can be transformed into the Painlev\'e V equation.
In (\ref{third order ode}), the $P_i$'s and $Q_1$ are
polynomials in x,
 \bean
 P_2(x)&:=&
  \left\{\begin{array}{l}
(4a_0a_2-a_1^2)x^2+2(2a_1a_2n-2a_0b_2+a_1b_1-2a_2b_0)x\\
  -4a_2^2n^2-4 (2a_1b_2-a_2b_1)n
  +4b_0b_2-b_1^2+2a_2^2\\
 \end{array}\right\}
\\ \\
   P_1(x)&:=&
    \left\{\begin{array}{l}
(2a_1a^2_2n+2a_0a_2b_2-a^2_1b_2+a_1a_2b_1-2a^2_2b_0)x\\
 -4a_2^3n^2+(-6a_1a_2b_2+4a_2^2b_1)n\\
 -2a_0b_2^2+a_1b_1b_2+2a_2b_0b_2-a_2b_1^2
\end{array}\right\}
 \\  \\
 Q_1(x)&:=&
  \left\{\begin{array}{l}
(2a_0a^2_2n-a_1^2a_2n+a_0a_1b_2-2a_0a_2b_1+a_1a_2b_0)x \\
+2a_1a^2_2n^2+(2a_1^2b_2-a_1a_2b_1-2a_2^2b_0)n
%    +(a_0b_1b_2-2a_1b_0b_2+a_2b_0b_1)
   \\
 +a_0b_1b_2-2a_1b_0b_2+a_2b_0b_1
\end{array}\right\}
 .\eean

\end{proposition}

The proof of Proposition 2.1 hinges on the following
Lemma, which we state in its full generality, although
only the case $I=E$ will be used; see \cite{AvM1, vM}.

\begin{lemma} Given a disjoint union of intervals
$I={\bigcup^r_1}[c_{2i-1},c_{2i}]\subset E$, the
integral
$$
\tau_n(t,c)=\int_{I^n}\Delta^2(z)\prod^n_1
e^{\sum_{i=1}^{\iy} t_iz_k^i}\rho(z_k)dz_k
$$
with $\rho(z)$ and $E$ as in (\ref{rational
logderivative}) and
 (\ref{boundary condition})
satisfies
\begin{description}

  \item[(i)] \underline{Virasoro constraints}
for all $m\geq -1$:%for $\beta =2$,
 \be
\left(-\sum_1^{2r} c_i^{m+1}A(c_i)\frac{\pl}{\pl
c_i}+{\mathbb V_m}
 %\sum_{i\geq 0}\left( a_i~ \BJ_{k+i,n}^{(2)}-b_i ~
% \BJ_{k+i+1,n}^{(1)}\right)
\right)\tau_n(t,c) =0.
 \ee
 with
$$
 {\mathbb V_m}
  :=\sum^2_{k=0}\left\{\begin{array}{c}
a_k(J^{(2)}_{k+m}+2n~J^{(1)}_{k+m}+n^2~J^{(0)}_{k+m})\\
 \\
 -b_k(J^{(1)}_{k+m+1}+n~\delta~J^{(0)}_{k+m+1})
\end{array}
\right\} ,
 $$
  where
 \bean
  J^{(2)}_k&=&\sum_{i+j=k}\frac{\pl^2}{\pl t_i\pl
t_j}+\sum_{-i+j=k}it_i
\frac{\pl}{\pl t_j}+\frac{1}{4}\sum_{-i-j=k}it_i ~jt_j\\
 \\
 J^{(1)}_k&=&\frac{\pl}{\pl t_k}+ \frac{1}{2}(-k)t_k,\quad
 J^{(0)}_k=\delta_{k0}.
 \eean

 \item[(ii)] The \underline{KP-hierarchy}\footnote{Given a polynomial
 $p(t_1,t_2,...)$, define the
customary Hirota symbol $p(\pl_t)f\circ g:=
p(\frac{\pl}{\pl y_1},\frac{\pl}{\pl
y_2},...)f(t+y)g(t-y) \Bigl|_{y=0}$. The ${\bf
s}_{\ell}$'s are the elementary Schur polynomials
$e^{\sum^{\iy}_{1}t_iz^i}:=\sum_{i\geq 0} {\bf
s}_i(t)z^i$ and for later use, set ${\bf
s}_{\ell}(\tilde \pl):={\bf s}_{\ell}(\frac{\pl}{\pl
t_1},\frac{1}{2}\frac{\pl}{\pl
t_2},\ldots).$}($k=0,1,2,\ldots$)
  $$
\left({\bf s} _{k+4}\bigl(\frac{\pl}{\pl
t_1},\frac{1}{2}\frac{\pl}{\pl
t_2},\frac{1}{3}\frac{\pl}{\pl
t_3},\ldots\bigr)-\frac{1}{2}\frac{\pl^2}{\pl t_1\pl
t_{k+3}}\right)\tau_n \circ\tau_n=0, $$
 of
which the first equation reads: \be \hspace{-1cm}
\left(\left(\frac{\pl}{\pl t_1}
\right)^4+3\left(\frac{\pl}{\pl
t_2}\right)^2-4\frac{\pl^2}{\pl t_1 \pl
t_3}\right)\log\tau_n+6\left(\frac{\pl^2}{\pl
t^2_1}\log\tau_n \right)^2=0. \label{KP-equation}\ee

\end{description}

\end{lemma}

\proof The most transparent way to prove this lemma is
via vector vertex operators, for which the
$\beta$-integrals
 \be
\tau_n(t,c;\beta):=\int_{E^n}|\Dt_n(x)|^{2\beta}\prod_{k=1}^n
\left(e^{\sum_1^{\iy}t_i x_k^i}\rho(x_k)dx_k\right), ~~
\mbox{for}~ n>0 \label{Integral with boundary}
 \ee
  are fixed points (see \cite{AvM1}).
 Another method, more computational, but much less conceptual, is to use
  a self-similarity argument, as in \cite{AvM1}. Namely, setting
 $$
 d\tau_n(z):=|\Dt_n(z)|^{2\beta}\prod_{k=1}^n
\left(e^{\sum_1^{\iy}t_i z_k^i}\rho(z_k)dz_k\right),$$
we have the following variational formula:
 \be
  \left.\frac{d}{d\vr}d\tau_n (z_i\mapsto
z_i+\vr A(z_i)z_i^{k+1} )\right|_{\vr=0}
=\sum^{\iy}_{\ell=0}
 \left(a_{\ell}~{}^{\beta}\!\BJ^{(2)}_{k+\ell,n}
-b _{\ell}~{}^{\beta}\!\BJ^{(1)}_{k+\ell
+1,n}\right)d\tau_n,
  \label{variational formula}\ee
with
\bean
%\lefteqn{ ~{}^{\beta}\!\BJ_{k,n}^{(2)}(t,n)
%}\\
 ~{}^{\beta}\!\BJ_{k,n}^{(2)}(t,n)
 &=& {\beta} ~~ {}^{\beta}\!\!J_k^{(2)} +
\Bigl(2n\beta +\!\!(k+1)\!(1\!-\! {\beta} )\Bigr) ~
{}^{\beta}\!\!J_k^{(1)} \!+ n\Bigl((n\!-\!1) {\beta}
+1\Bigr) \dt_{k0},\\
{}^{\beta}\!\BJ_{k,n}^{(1)}(t,n)&=& ~
{}^{\beta}\!\!J_k^{(1)}+n\dt_{k0},
% ~\mbox{with}~ ~ {}^{\beta}\BJ_{k,n}^{(0)}=nJ_k^{(0)}
% = n\dt_{0k}
% \\
% {}^{\beta}\!\BJ_{k,n}^{(0)}&=&nJ_k^{(0)}= n\dt_{0k}
%
 \eean
% $$ ~
%{}^{\beta}\!\BJ_{k,n}^{(1)}(t,n)= ~
%{}^{\beta}\!\!J_k^{(1)}+nJ_k^{(0)}
% ~\mbox{and}~ ~ {}^{\beta}\!\BJ_{k,n}^{(0)}=nJ_k^{(0)}= n\dt_{0k}
%$$
  where
  \bea
%  J_k^{(0)}&=& \dt_{k0}\\
 ~ {}^{\beta}\!\! J_k^{(1)}&=&\frac{\pl}{\pl
t_k}+\frac{1}{2\beta}(-k)t_{-k}
 \nonumber
 \\
  ^{\beta}\!\!J^{(2)}_{k}&=&\sum_{i+j=k}\frac{\pl^2}{\pl
 t_{i}\pl t_{j}}+\frac{1}{\beta}\sum_{-i+j=k}it_{i}\frac{\pl}{\pl
 t_{j}}+\frac{1}{4\beta^2}\sum_{-i-j=k}it_{i}jt_{j}.
%\nonumber
 \eea
 The change of
integration variable $ z_i\mapsto z_i+\vr
A(z_i)z_i^{k+1}
  $
   in the
integral (\ref{Integral with boundary}) leaves the
integral invariant, but it induces a change of limits of
integration, given by the inverse of the map above;
namely the $c_i$'s in
$I={\bigcup^r_1}[c_{2i-1},c_{2i}]$, get mapped as
follows $$ c_i\mapsto c_i-\vr A(c_i)c_i^{k+1}+O(\vr^2).
$$ Therefore, setting
$$I^{\vr}=\displaystyle{\bigcup^r_1}[c_{2i-1}-\vr
A(c_{2i-1}) c_{2i-1}^{k+1}+O(\vr^2),c_{2i}-\vr
A(c_{2i})c_{2i}^{k+1} +O(\vr^2)],$$ we find, using
(\ref{variational formula}) and the fundamental theorem
of calculus,
\begin{eqnarray*}
0\!&=&\frac{\pl}{\pl\vr}\int_{(I^{\vr})^{2n}}|
 \Delta_{2n}(z\!+\!\vr \! A(z)z^{k+1})|^{2\beta}
  \prod^{ n}_{i=1}e^{-V(z_i+\vr A(z_i)z_i^{k+1} ,t)}
   d(z_i\!+\!\vr\! A(z_i)z_i^{k+1})\\
&=&\left(-\sum^{2r}_{i=1}c_i^{k+1} A(c_i)\frac{\pl}{\pl
c_i} +\sum^{\iy}_{\ell=0}\left( a_{\ell}
~{}^{\beta}\!\BJ^{(2)}_{k+\ell,2n}
 -b_{\ell}~{}^{\beta}\!\BJ^{(1)}_{k+\ell+1,2n}
\right)\right)  \tau_{n}(t,c,\beta).
\end{eqnarray*}
For Lemma 2.2 one sets $\beta =1$; also, when $I=E$, the
condition (\ref{boundary condition}) implies that the
boundary terms in the formula above are absent. For
statement (ii), concerning the KP equation, we refer the
reader to \cite{vM}. This ends the proof of Lemma
2.2.\qed

{\medskip\noindent{\it Proof of Proposition 2.1:\/} }
 Setting $F(t):= F_n(t)=\log \tau_n(t)$, a few of the Virasoro constraints
 of Lemma 2.2, evaluated along the locus
 $$
 \LR:= \{t=(x,0,0,\ldots) \},
 $$
read as follows:
\bean
 \frac{\mathbb
V_{-1}\tau_n}{\tau_n}\Bigr\vert_{\LR}&=&a_0
\left(\sum_{i\geq 2}it_i\frac{\pl F}{\pl t_{i-1}}+ n
t_1\right)+a_1
\left(\sum_{i\geq 1}it_i\frac{\pl F}{\pl t_{i}}+ n^2\right)\\
& &~~~~+a_2\left(\sum_{i\geq 1}it_i\frac{\pl F}{\pl
t_{i+1}}+ 2n~\frac{\pl F}{\pl
t_1}\right)-\left(b_0n+b_1\frac{\pl F}{\pl t_1}+
b_2\frac{\pl F}{\pl t_2}\right) \Bigr\vert_{\LR}
 \\
&=&n(a_0x+a_1n-b_0)+(a_1x\!+\!2n a_2\!-\!b_1)\frac{\pl
F}{\pl
t_1}%\Bigr\vert_{\LR}
 +(a_2x-b_2)\frac{\pl F}{\pl
t_2}\Bigr\vert_{\LR}=0. \\&&\\&&\\
%  \eean
%
%
%
%\bigbreak
%\bean
%
\frac{\mathbb
V_{0}\tau_n}{\tau_n}\Bigr\vert_{\LR}&=&a_0
\left(\sum_{i\geq 1}it_i\frac{\pl F}{\pl t_{i}}+
n^2\right)+a_1 \left(\sum_{i\geq 1}it_i\frac{\pl F}{\pl
t_{i+1}}+ 2n\frac{\pl F}{\pl
t_1}\right)\\
& &~~~~+a_2\left(\sum_{i\geq 1}it_i\frac{\pl F}{\pl t_{i+2}}+
\frac{\pl^2 F}{\pl t^2_1}+\left(\frac{\pl F}{\pl
t_1}\right)^2+2n\frac{\pl
F}{\pl t_2}\right)\\
& &~~~~-\left(b_0\frac{\pl F}{\pl t_1}
 +b_1\frac{\pl F}{\pl t_2}+
b_2\frac{\pl F}{\pl t_3}\right)\Bigr\vert_{\LR}%\\
 \eean\bean
&=&a_0n^2+(a_0x+2n a_1-b_0)\frac{\pl F}{\pl
t_1}+a_2\left(\frac{\pl F}{\pl t_1}\right)^2+ a_2 \frac{\pl^2 F}{\pl t_1^2}\\
& &~~~~+(a_1x+2n a_2-b_1)\frac{\pl F}{\pl
t_{2}}+(a_2x-b_2) \frac{\pl F}{\pl
t_3}\Bigr\vert_{\LR}=0. \\&&\\&&\\ %\eean
%
%
%
%\bigbreak \bean
%
\frac{\pl}{\pl t_1}\frac{\mathbb
V_{-1}\tau_n}{\tau_n}\Bigr\vert_{\LR}&=&a_0n+ \sum^2_1
a_i\frac{\pl F}{\pl t_{i}}+ (a_1x+2n a_2-b_1)
\frac{\pl^2 F}{\pl t_{1}^2}\\
& &~~~~+(a_2x-b_2) \frac{\pl^2 F}{\pl t_1\pl
t_2}\Bigr\vert_{\LR}=0. \eean

%\bigbreak

 \bean
\frac{\pl}{\pl t_1}\frac{\mathbb
V_{0}\tau_n}{\tau_n}\Bigr\vert_{\LR}
 &=&
  \sum^2_{i=0}  a_i\frac{\pl F}{\pl t_{i+1}}
   + a_2\left(\frac{\pl^3 F}{\pl t_1^3}
 +2\frac{\pl F}{\pl t_1}\frac{\pl^2 F}{\pl t_1^2}\right)\\
& &~~~~+(a_0x+2n a_1-b_0) \frac{\pl^2 F}{\pl
t_1^2}+(a_1x+2n-b_1)\frac{\pl^2F}{\pl t_1\pl t_2}\\
& &~~~~~~~~+(a_2x-b_2) \frac{\pl^2 F}{\pl t_1\pl
t_3}\Bigr\vert_{\LR}=0. \\&&\\&&\\ %\eean
%
%
%
%\bigbreak \bean
%
\frac{\pl}{\pl t_2}\frac{\mathbb
V_{-1}\tau_n}{\tau_n}\Bigr\vert_{\LR}&=&2 \sum^2_{i=0}
a_i\frac{\pl F}{\pl t_{i+1}}+ (a_1x+2n a_2-b_1)
\frac{\pl^2 F}{\pl t_{1}\pl t_2}\\
& &~~~~+(a_2x-b_2) \frac{\pl^2 F}{\pl
t_2^2}\Bigr\vert_{\LR}=0. \eean

\bigbreak

\noindent These five equations form a linear system in
the five unknowns
$$
\frac{\pl F}{\pl t_{2}}\Bigr\vert_{\LR}, \frac{\pl^2
F}{\pl t_{1}\pl t_2}\Bigr\vert_{\LR},\frac{\pl^2 F} {\pl
t_{2}^2}\Bigr\vert_{\LR}, \frac{\pl F}{\pl
t_{3}}\Bigr\vert_{\LR},\frac{\pl^2 F}{\pl t_{1}\pl
t_3}\Bigr\vert_{\LR} ,
 $$
 which upon solving in terms of $\left(
  \frac{\pl  }{\pl t_{1}} \right)^k  F  \Bigr\vert_{\LR}$
  and substituting into the KP-equation (remembering
  equation (\ref{KP-equation}))
$$
\left( \frac{\pl^4 }{\pl t_{1}^4} +3 \frac{\pl^2 }{\pl
t_{2}^2}  -4\frac{\pl^2}{\pl t_{1}\pl
t_3}\right)F+6\left(\frac{\pl^2 F} {\pl
t_{1}^2}\right)^2=0
$$
yields a differential equation in $F$. Pure $F$
 never appears in this equation, because the Virasoro
constraints only contain partials of $F$. Therefore, it
is a differential equation in $g(x)= \frac{\pl}{\pl t_1}
F(t_1,0,0,\ldots)\Bigr|_{t_1=x}$, which one computes has
the form (\ref{third order ode}); this ends the proof of
Proposition 2.1. \qed

The proof of Theorem 1.5 will be given at the end of
section 3.

%\newpage

\section{Hermitian matrix integrals and Painlev\'e equations}

Define
 $$ \HR_n=\{n\times n\mbox{~Hermitian
matrices}\},
$$
and, for an interval $I\subset \BR$,
$$ \HR_n(I)=\{M\in \HR_n \mbox{~with spectrum in $I$}\}.
$$
We give a proof of the following theorem due to
Forrester and Witte \cite{Forrester Witte1,Forrester
Witte2}, using Proposition 2.1.

\begin{theorem}
 For\footnote{For the first integral of (\ref{matrix
integral}), one can choose the intervals
$I_1=(-\iy,s],[s,\iy)$ or $(-\iy,\iy)$. For the second
integral of (\ref{matrix integral}), one may choose the
intervals $I_2=[0,s],[s,\iy)$ or $[0,\iy)$
  }
   $a,b>-1$, the logarithmic derivatives of the integrals
 \bea
% E_1(c)
 h(s)&:=& \frac{d}{ds}\log\int_{\HR_n(-\iy,s]
       }
  \det(M-sI)^a e^{-\Tr M^2} dM
  \no\\
%   E_2(c)
   k(s)&:=& s\frac{d}{ds}\log \int_{\HR_n[0,s]
       }
  \det(sI-M)^b \det M^a  e^{-\Tr M} dM
 \label{matrix integral}\eea
% can be expressed as %satisfy
%  \bea
%  E_1(c)&=&E_1(0)\exp{ \int_0^c h(y) dy} \no\\  \no\\
%  E_2(c)&=&E_2(1)\exp{ \int_1^c \frac{k(y)}{y} dy}
%  ,\eea
%  with $h(y)$ and $k(y)$
satisfy the Painlev\'e IV and V
  equations, respectively:
 \be
h^{\prime\prime\prime}+6h^{\prime
2}-4(s^2+2(a-n))h^{\prime}+4sh-8an=0 \label{ode in h}
 ,\ee
  and
 \bea
 \lefteqn{k^{\prime\prime\prime}+\frac{k^{\prime\prime}}{s}
  +\frac{6}{s}k^{\prime 2}
 -\frac{4}{s^2}kk'-\left(s^2-2s(2n+a-b)+ (a+b)^2\right)
  \frac{k'}{s^2}}   \no \\
& &~~~~~~~~~~~~~~~~
 +(s-2n-a+b)\frac{k}{s^2}-\frac{bn}{s^2}(s+a+b)=0.
\label{ode in k}\eea

\end{theorem}

{\medskip\noindent{\it Proof of Theorem 3.1:\/} }
Set
 \be
 \tau(x)=\int_{I^n} \Dt^2(z) \prod_1^n e^{xz_k}\rho(z_k)dz_k
 ,\label{tau in x}\ee
  with $\rho$ and $I$ as in
  (\ref{rational logderivative})
 and (\ref{boundary condition}).

 {\bf (i)} Then, expressed in spectral coordinates and
 making the substitution $y_i=z_i-s$ for $1\leq
 i\leq n$, the first matrix integral in (\ref{matrix
 integral}) reads, for $I=(-\iy,s],[s,\iy)$ and
 $(-\iy,\iy)$,
 \bea
 %E_1(s) &=&
 \int_{(I_1)^{n}}
  \Dt^2(z) \prod_1^n (z_i-s)^a e^{-z_i^2}dz_i
  %\no\\
        &=& e^{-ns^2} \int_{(I'_1)^n}
         \Dt^2(y) \prod_1^n y_i^ae^{-y_i^2-2sy_i}dy_i
         \no\\
         &=& e^{-ns^2}\tau(-2s)
   \eea
with $I'_1=(-\iy,0],[0,\iy),(-\iy,\iy)$. Here
  $\tau(x)$ as in (\ref{tau in x}) contains
$\rho(z)=z^ae^{-z^2}$, for which
$$
-\frac{\rho'}{\rho}=\frac{-a+2z^2}{z}
 .%~~\mbox{and}~I=(-\iy,0],(-\iy,\iy)~\mbox{or}~[0,\iy).
$$
 Thus, setting
  \bea
& &a_0=a_2=0,\quad a_1=1\nonumber\\
& &b_0=-a,\quad b_1=0,\quad b_2=2.
 \label{special values of a}\eea
%
%Setting (\ref{special values of a})
 in the equation
(\ref{third order ode}), one deduces that
% for
% $I=(-\iy,0]$,$~(-\iy,\iy)$, $[0,\iy)$,
%
$$
g(x):= \frac{\pl}{\pl x}\log \tau(x)
 %\int_{I'^n}\Dt^2(z)
%\prod_1^n \left( e^{-z_k^2+xz_k}z_k^adz_k\right), ~
 ,~~\mbox{for}~
I =\left\{  \begin{array}{l} (-\iy,0]\\ (-\iy,\iy) \\ (0,\iy)\\
  \end{array}
   \right.
   , $$
%
%
%$$
%g=\frac{\pl F}{\pl x}
%$$
%
satisfies
$$
g^{\prime\prime\prime}+6g^{\prime
2}-g^{\prime}\left(\frac{x^2}{4}+4n+2a\right)+\frac{xg}{4}+\frac{n}{2}(n+a)=0.
$$
 But we need a differential equation for
 $e^{-ns^2}\tau(-2s)$, instead of $\tau(x)$. Therefore
 consider
 $$
h(s)=\frac{\pl}{\pl
s}\log(e^{-ns^2}\tau(-2s))=-2ns-2g(x)\Bigr\vert_{x=-2s}
, $$
which relates to $g(x)$ as follows,
\bean
h'(s)&=&-2n+4g'(x)\Bigr\vert_{x=-2s}\\
h''(s)&=&-8g''(x)\Bigr\vert_{x=-2s}\\
h'''(s)&=&16g'''(x)\Bigr\vert_{x=-2s}
 .\eean
%{ode in h}
Expressing $g(x),g'(x),g''(x),g'''(x)\Bigr\vert_{x=-2s}$
in terms of $h,h',h''$ and $h'''$, and setting $x=-2s$
yield the differential equation (\ref{ode in h}), which
according to the table in Appendix 1 is a version of
Painlev\'e IV; this establishes the first part of
Theorem 3.1.
%
%$$
%h^{\prime\prime\prime}+6h^{\prime
%2}-4h^{\prime}(c^2+2(a-n))+4ch-8an=0.
%$$

\bigbreak

{\bf (ii)} Then, making the substitution $y_i=z_i/s$
 for $1\leq i\leq n$,
 \bea
% E_2(s) &=&
\lefteqn{\int_{(I_2)^{n}}
  \Dt^2(z) \prod_1^n (s-z_i)^b z_i^a e^{-z_i}dz_i
 } \no\\
        &=&  s^{n(n+a+b)} \int_{(I'_2)^n}
         \Dt^2(y) \prod_1^n (1-y_i)^b y_i^ae^{-sy_i }dy_i
         \no\\
         &=& s^{n(n+a+b)}\tau(-s)
   \eea
with $I_2=[0,s], ~I'_2=[0,1]$, for $I_2=[s,\iy],
~I'_2=[1,\iy]$ and finally for $I_2=[0,\iy],
~I'_2=[0,\iy]$;
 $\tau(x)$ now corresponds to
$\rho(z)=z^a(1-z)^b$, and so
$$
-\frac{\rho'}{\rho}=-\frac{a}{z}+\frac{b}{1-z}=\frac{a-(a+b)z}{-z+z^2}
.$$
Thus
\bean
\begin{array}{lll}
a_0=0,&a_1=-1,&a_2=1\\
 \\
b_0=a,&b_1=-(a+b),&b_2=0.
\end{array}
\eean
Setting these special values in the equation (\ref{third
order ode}), one checks that
% for
% $I=(-\iy,0]$,$~(-\iy,\iy)$, $[0,\iy)$,
%
$$
g(x):= \frac{\pl}{\pl x}\log\int_{I_2^{\prime
n}}\Dt^2(z) \prod_1^n
  e^{xz_k}z_k^a  (1-z_k)^b dz_k  , ~
% \mbox{for}~
%I=\left\{  \begin{array}{l} (-\iy,0]??\\ (-\iy,\iy) \\ (0,\iy)\\
%  \end{array}
%   \right.
    $$
%
%
%$$
%g=\frac{\pl F}{\pl x}
%$$
%
satisfies
\bea  \lefteqn{g^{\prime\prime\prime}+6g^{\prime
2}+\frac{4}{x}(g''+2gg')+2\left(\frac{g}{x}\right)^2}\no\\
 \no\\
&-&\frac{g'}{x^2}\left(x^2+2(2n+a-b)x+(2n+a+b)^2-2\right)
\no\\
 \no\\
&-&\frac{g}{x^3}\left((2n+a-b)x+(2n+a+b)^2\right)\no\\
 \no\\
&+&\frac{n}{x^3}(n+a)(x+2n+a+b)=0. \label{ode2}\eea

\bigbreak

Since $g(x)=\displaystyle{\frac{\pl}{\pl
x}}\log\tau(x)$, we have
\bean
k(s)&=&s~\frac{\pl}{\pl s}\log s^{n(n+a+b)}\tau(-s)\\
&=&n(n+a+b)-sg(-s) \eean and so

\bean
g(-s)&=&\frac{1}{s}\left( {-k(s)+n(n+a+b)} \right)\\
 \\
g'(-s)&=&\frac{1}{s^2}\left(sk'(s)-k(s)+n(n+a+b)\right)\\
 \\
g''(-s)&=&\frac{2}{s^3}\left(-s^2k''(s)/2+sk'(s)-k(s)
    +n(n+a+b)\right)\\
 \\
g'''(-s)&=&\frac{3}{s^4}\left(s^3k'''(s)/3-s^2k''(s)
 +2sk'(s)-2k(s)+2n(n+a+b)
\right).\\
\eean
Substituting these expressions into (\ref{ode2}) yields
 the differential equation (\ref{ode in k}), which again
  referring to the table in appendix 1 is Painlev\'e V,
  ending the proof of Theorem 3.1.
  %Note this also proves Theorem 1.5.
 \qed

{\medskip\noindent{\it Proof of Theorem 1.5:\/} }
Because of the form (\ref{rational logderivative}) of
$\rho(z)dz$, we have for an open set of real constants
 $a_i$ and $b_i$,
$$
 \rho(z)=   (z-a)^{\al} (b-z)^{\beta}e^{\gamma z}
 ,\mbox{  with  } \al, \beta >-1  ~~\mbox{and}~~a,b \in \BR,$$
  which, upon making a linear change of variables $z\mapsto y$
   in
 the integral (\ref{integral g}), leads to
 $$
 g(x) =na+\frac{\pl}{\pl
x}\log\int_{[0,1]^n}\Dt^2(y) \prod_1^n e^{ (x+\gamma
)(b-a)y_k}y_k^{\al } (1-y_k)^{\beta }  dy_k
$$
and so
 $$
 \frac{1}{b-a}\left(g\Bigl(\frac{x'}{b-a}-\gamma\Bigr)-na\right)
  =  \frac{\pl}{\pl
x'}\log\int_{[0,1]^n}\Dt^2(y) \prod_1^n e^{
x'y_k}y_k^{\al } (1-y_k)^{\beta }  dy_k ,
 $$
  which is
shown to be a solution of Painlev\'e V in part {\bf
(ii)} of the proof of Theorem 3.1 in this section.

To compute the constant %$\dt$
in (\ref{1.0.21}) (that is $\dt$ in equation
(\ref{Cosgrove})), multiply the equation (\ref{1.0.21})
with $P^2/4$, set $y=0$ and use the explicit expressions
for the polynomials $P,~Q,~R$ and the function $f(y)$,
as in (\ref{f(y)}), yielding the identity

\be\dt= -\left[ Q f'+R\right] f' + \left[f'^2+Q' f' +R'
\right] f  - \frac{1}{2}Q'' f ^2 \Big
|_{y=0}.\label{a}\ee
 All the
expressions can readily be computed, except for $f'(0)$,
which requires some argument. One computes
\bean
f'(0)&=&
 \frac{1}{b-a}\left(na-g(-\gamma)\right)\\
&=&\frac{n}{b-a} \left(  a-
  \frac{\int_{[a,b]^n}\Dt^2(z)(\frac{1}{n}\sum_1^n z_k)
  %\prod_1^n e^{xz_k}
  (z_k-a)^{\al} (b-z_k)^{\beta}dz_k}
{ \int_{[a,b]^n}\Dt^2(z) \prod_1^n %e^{xz_k}
  (z_k-a)^{\al}
(b-z_k)^{\beta}dz_k} \right)
\\
&=&
 \frac{n}{b-a} \left(  a- \la z_1 \ra_{_{(a,b)}}\right)
\\
&=& -\frac{n(n+\al)}{2n+\al +\beta}. \eean
To establish the last equality above, we need the Aomoto
extension \cite{A} (see \cite{AvM2}, Appendix D) of
Selberg's integral:%\footnote{where $Re\,\alpha$,
%$Re\,\beta >-1$, $Re\,\gamma >
%-\min\displaystyle{\left(\frac{1}{n}, \frac{Re\,\alpha
%+1}{n-1}, \frac{Re\,\beta +1}{n-1}\right)}$}

\bea
 \la x_1...x_m\ra_{_{(0,1)}}&:=&
 \displaystyle{\frac{\displaystyle{\int_{[0,1]^n}}
 x_1\ldots
x_m\left|\Delta(x)\right|^{2\gamma}\prod^n_{j=1}x^{\alpha}_{j}
(1-x_{j})^{\beta}dx_j}{\displaystyle{\int_{[0,1]^n}}\left|\Delta(x)\right|^{2\gamma}
\prod^n_{j=1}x^{\alpha}_{j}
(1-x_{j})^{\beta}dx_j}}\nonumber\\
 &=&\prod^m_{j=1}\frac{\alpha
+1+(n-j)\gamma}{\alpha+\beta+2+(2n-j-1)\gamma}. \eea In
particular, setting $\gamma=1$, this formula implies
 \be \la
z_1\ra_{_{(a,b)}}=\la a+(b-a)x_1 \ra_{_{(a,b)}} =a+
(b-a)\la x_1 \ra_{_{(0,1)}} =a+
(b-a){{n+\alpha}\over{2\,n+\beta+\alpha}} ~~ \ee
Adding up all the pieces in (\ref{a}), one finds the
value of the constant $\dt$ in equation
(\ref{Cosgrove})%(\ref{1.0.21})
, ending the proof of
Theorem 1.5.\qed

%\newpage

%\section{Probability on words and partitions}

%

%\newpage

\section{A limit theorem for the probability on words and Chi-square}

The purpose of this section is to prove Theorem 1.2 and
Corollary 1.3. Given $m$ independent normal
$N(0,1)$-random variables $X_i$, the sum
$$
Z_m=\sum_1^mX_i^2
$$
is $\chi^2_m$-distributed, with Fourier transform
$$
E(e^{t\,Z_m})=(1-2t)^{-m/2}.
$$
Developing both sides in $t$ yields the moments
$$
\frac{1}{k!}E\left(\left(\frac{1}{2}Z_m\right)^k\right)=
 \left(\begin{array}{c}\frac{m}{2}-1+k\\k\end{array}\right).
$$

{\medskip\noindent{\it Proof of Theorem 1.2:\/} }
 From Stirling formula $N!=\sqrt{2\pi N}~
e^{N(\log N -1)}(1+{\bf O}(\frac{1}{N}))$ and
$(N+\alpha)!\simeq N!N^{\alpha}$, we have for large $n$,
\bea
\begin{array}{cll}
(pn+k-p^2)!&\simeq (pn)!(pn)^{k-p^2}&
 \simeq
  \displaystyle{\frac{(n! ~ p^n)^p}{(2\pi n)^{\frac{p-1}{2}}}\sqrt{p}
 ~ (pn)^{k-p^2}}
 \no \\
  \no\\
\displaystyle{\prod^p_1}  (n-j)!   &\simeq
  \displaystyle{\prod^p_1}  (n!~ n^{-j}) &\simeq
 (n!)^p ~ n^{-\frac{p(p+1)}{2}}
 \end{array}
  \eea  \vspace{-1 cm}\be \label{stirling}\ee
  meaning here that the ratios between two consecutive expressions
  tend to $1$, when $n\rg
  \iy$.
On the one hand, we have, changing the summation index
$\ell \rg \ell -p(n-p)$,
\bea
\lefteqn{\lim_{n\rg\iy}\prod^p_{j=1}\frac{(n-j)!}{(q-j)!}
\sum_{\ell\geq p(n-p)}\frac{p^{\ell}(nu)^{\ell
-p(n-p)}}{\ell !}E^{\ell,p}\left(
I_{\lb\supseteq\mu}(\lb)\prod_{n-q<j\leq
n-p}h^{\lb}_{ij}\right) }
\no\\\no\\
&=&\lim_{n\rg\iy}\prod^p_{j=1}\frac{(n-j)!}{(q-j)!}\sum_{k=0}^{\iy}
\frac{p^{p(n-p)}(pnu)^k}{(k+p(n-p))!}E^{k+p(n-p),p}
\left( I_{(\lb\supseteq\mu)}(\lb)
 \prod_{n-q<j\leq n-p}h^{\lb}_{ij}\right)
\no\eea\bea% \no\\ \no\\
&=&\!\!
\sum_{k=0}^{\iy}u^k\!\!\lim_{n\rg\iy}\frac{(pn)^kp^{p(n-p)}}
 {\displaystyle{\prod^p_1}(q-j)!}\frac{\displaystyle{
\prod^p_1}(n-j)!}{(pn+k-p^2)!}
E^{k+p(n-p),p}\left( I_{\lb\supseteq\mu}(\lb)
\prod_{n-q<j\leq n-p}h^{\lb}_{ij}\right)\no\\
\no\\
&=&\sum_{k=0}^{\iy}u^k\lim_{n\rg\iy}
\frac{(2\pi)^{\frac{p-1}{2}}n^{\frac{p^2-1}{2}}}
{\displaystyle{\prod^p_1}(q-j)!\sqrt{p}}~~
E^{k+p(n-p),p}\left(
I_{\lb\supseteq\mu}(\lb)\prod_{n-q<j\leq
n-p}h^{\lb}_{ij}\right),
 \label{expression 1}\eea
  using (\ref{stirling}) in the last equality.

On the other hand, using (\ref{restiction of
hypergeometric}), (\ref{symbol}), (\ref{proposition4.1})
and
\bea \lim_{n\rightarrow \iy}\frac{(n)_{\lb}}{n^{|\lb|}}
  &=&
  \lim_{n\rightarrow \iy}\left[
  \frac{n(n+1)\ldots(n+\lb_1-1)}{n^{\lb_1}}\right]
  \left[
  \frac{(n-1)n \ldots(n+\lb_2-2)}{n^{\lb_2}}\right]\ldots
,%
 \no\\
 &=& 1
 \eea
 we have
 \bea
\lefteqn{\hspace{-3cm}\lim_{n\rg\iy}\prod^p_{j=1}\frac{(n-j)!}{(q-j)!}\sum_{\ell\geq
p(n-p)}\frac{p^{\ell}(nu)^{\ell -p(n-p)}}{\ell
!}E^{\ell,p}\left(I_{\lb\supseteq\mu}(\lb)\prod_{n-q<j\leq
n-p}h^{\lb}_{ij}\right)}\no\\
\no\\
&=&
 \lefteqn{\lim_{n\rg\iy}\,_2F_1^{(1)}(p,q;n;x)
 \Big\vert_{\sum_{\ell}x_{\ell}^i=nu\,\dt_{i1}}}\no\\
\no\\
&=&\lim_{n\rg\iy}\sum_{k=0}^{\iy}u^k\sum_{\kappa\in\BY_k}
 \frac{n^k(p)_{\kappa}(q)_{\kappa}}
 {(h^{\kappa})^2(n)_{\kappa}}\no\\
\no\\
&=&\sum_{k=0}^{\iy}u^k\sum_{\kappa\in\BY_k}
 \frac{(p)_{\kappa}(q)_{\kappa}}{(h^{\kappa})^2}\no\\
\no\\
&=&\sum_{k=0}^{\iy}u^k\sum_{\kappa\in\BY_k}
 {\bf s}_{\kappa}(1^p)
 {\bf s}_{\kappa}(1^q)
 ,~~\mbox{using (\ref{semi-standard})}
  \no\\
  &=&  \sum_{\kappa\in\BY }
 {\bf s}_{\kappa}(u^p)
 {\bf s}_{\kappa}(1^q)
\no\\
%&=&\,_2F_0^{(1)}(p,q;x)
%\Big\vert_{\sum_{\ell}x_{\ell}^i=u\,\dt_{i1}}\no\\
\no\\
&=&(1-u)^{-pq}, \mbox{  using the Cauchy identity}\no\\
\no\\
&=&E\left(e^{\frac{1}{2}u\,Z_{2pq}}\right)
 \label{expression 2}\eea
Then comparing the coefficients of $u^k$ in the
identical asymptotic expressions (\ref{expression 1})
and (\ref{expression 2}) yields Theorem 1.2. \qed

{\medskip\noindent{\it Proof of Corollary 1.3:\/} }
Setting $p=q$ and $N= n-p $ in Theorem 1.2, the
expectation in (\ref{Theorem 5.1}) becomes
$P^{Np+k,p}(\lb \supseteq \mu)$ for fixed $\mu$. By RSK
and (\ref{RSK}), the condition $\lb \supseteq \mu$
translates into $d_1(\pi)=\lb_1^{\top}=p$ and
$i_{p-1}(\pi)= \sum_1^p\lb_i-\lb_p=Np+k-\lb_p$, with
$\lb_p\geq N$. This means $i_{p-1}(\pi)\leq
k+(p-1)N$.\qed

%\newpage

\section{A limit theorem for the Poissonized probability on words}

In the next Theorem, we consider the expectation for the
Poissonized probability $P_{x,p}$ on Young diagrams
$\lb\in\BY_{\ell}^{(p)}$, as defined in
(\ref{Poissonized probability}). Before proving Theorem
1.4, one needs the following proposition:

\begin{proposition}   (\cite{Johansson1,TW})%, theorem 1.6)}
 For every continuous function $g$
on $\BR^p$, we have

\bean \lefteqn{\lim_{\ell\rg\iy}
E^{\ell,p}\left(g\left(\frac{\lb_1-\ell/p}{\sqrt{2\ell/p}},\ldots,
\frac{\lb_p-\ell/p}{\sqrt{2\ell/p}}\right)\right)\qquad\vert
\lb\vert =\ell} \\
& &  \\
&=&p!\sqrt{\pi p}\renewcommand{\arraystretch}{0.5}
\begin{array}[t]{c}
\int\\
{\scriptstyle x\in\BR^p}\\
{\scriptstyle x_1>\ldots>x_p}\\
{\scriptstyle \sum_1^p x_i=0}
\end{array}
\renewcommand{\arraystretch}{1}g(x_1,\ldots,x_p)\vp_p(x_1,\ldots,x_p)dx_1\ldots
dx_{p-1} \eean
%
%\bigbreak
and
 \bean \lefteqn{\lim_{x\rg\iy}
E_{x,p}\left(g\left(\frac{\lb_1-x}{\sqrt{2x}},\ldots,
\frac{\lb_p-x}{\sqrt{2x}}\right)\right)} \\
& &  \\
&=&p!\renewcommand{\arraystretch}{0.5}
\begin{array}[t]{c}
\int\\
{\scriptstyle x\in\BR^p}\\
{\scriptstyle x_1>\ldots>x_p}
\end{array}
\renewcommand{\arraystretch}{1}
g(x_1,\ldots,x_p) \vp_p(x_1,\ldots,x_p)dx_1\ldots dx_p
 ~,\eean
  where $\vp_p(x_1,\ldots,x_p)$ is the probability
  density
$$
\vp_p(x_1,\ldots,x_p)=\frac{1}{Z_p}\Delta_p(x)^2\prod^p_{j=1}e^{-x^2_j},
$$
with
$$
Z_p=(2\pi)^{p/2}~~\frac{2^{-p^2/2}}
{\displaystyle{\prod^p_1}j!}.
 $$
\end{proposition}

\noindent{\it Proof of Theorem 1.4:\/}   Given partition
$\lb$ with the probability measure (\ref{Poissonized
probability}), consider the random variable

 \be
\vr_i(\lb):=\frac{\lb_i-x}{\sqrt{2x}}.
\label{epsilon}
 \ee
  From (\ref{rescaling}), we have
\be
 s=\frac{(n-p)-x}{\sqrt{2x}}
 .\label{c}\ee
Consider the hook length $h^{\lb}_{(ij)}$ for box $(i,j)
\in \lb$, such that $n-q<j\leq n-p$. Define $r_{ij}$
such that
\bean
h^{\lb}_{(ij)}&:=&\lb_i-(n-p)+r_{ij}\\
&=&\sqrt{2x}(\vr_i-s)+r_{ij}, \eean
 upon using formulas (\ref{epsilon}) and (\ref{c})
  in the last equality. The position $(i,j)$ in the
  partition $\lb$ such that $n-q<j\leq n-p$ implies that
  $r_{ij}\leq (q-p)+p=q$, from visual inspection of Figure 1.
  We also have
\bean
\{\lb\in\BY^{(p)},\lb\supseteq\mu\}&=&\{\lb\in\BY^{(p)},\lb_i\geq
n-p,\mbox{all~}1\leq i\leq p\}\\
& &\\
&=&\left\{\lb\in\BY^{(p)},\frac{\lb_i-x}{\sqrt{2x}}\geq\frac{(n-p)-x}{\sqrt{2x}}
,\mbox{all~}1\leq i\leq
 p\right\}\\
 & & \\
 &=&\bigcap^p_{i=1}\{\lb\in\BY^{(p)},\vr_i(\lb)\geq s\}.
\eean

 On the set
$\vr_i(\lb)\geq s$ and for $\lb\in\BY^{(p)}$,  we have,
since $0\leq r_{ij}\leq q$,
 \bean
\lefteqn{(2x)^{-\frac{p(q-p)}{2}}\prod_{(i,j)
 \in\lb\atop{n-q<j\leq
n-p}}h^{\lb}_{(i,j)}- \prod_{1\leq i\leq p} (\vr _
i-s)^{q-p}}\\
&=&
(2x)^{-\frac{p(q-p)}{2}}\prod_{(i,j)\in\lb\atop{n-q<j\leq
n-p}}(\sqrt{2x}(\vr_i-s)+r_{ij})- \prod_{1\leq i\leq p}
(\vr _ i-s)^{q-p}
% - \prod_{1\leq i\leq p} (\vr _ i-c)^{q-p}
 \eean\bean%\\
 &=&{ \prod_{1\leq i\leq p\atop{n-q<j\leq
 n-p}}
\left((\vr _ i-s)+\frac{r_{ij}}{\sqrt{2x}}\right)
  - \prod_{1\leq i\leq p} (\vr _ i-s)^{q-p}
}\\
&\leq&
 \prod_{1\leq i\leq p}
\left((\vr _ i-s)+\frac{q}{\sqrt{2x}}\right)^{q-p}
  - \prod_{1\leq i\leq p} (\vr _ i-s)^{q-p}
\\
 &=&
 \sum_{{\footnotesize\begin{array}{l}
       0\leq\ell_i\leq q-p \\
       1\leq i\leq p \\
       \sum\ell_i\geq 1
       \end{array}}}
    \left(\frac{q}{\sqrt{2x}}\right)^{\sum\ell_i}
    \prod_{i=1}^{p}\left({q-p}\atop {\ell_i}\right)
    \prod_{i=1}^p(\vr_i-s)^{q-p-\ell_i}
 \eean
The positive expression is now estimated as follows:
 \bean
 \lefteqn{
 (2x)^{-\frac{p(q-p)}{2}} E^{\ell,p}\left(I_{\{\lb\supseteq\mu\}}(\lb)
\prod_{(i,j)\in\lb\atop{n-q<j\leq
n-p}}h^{\lb}_{(i,j)}\right)
 }\\
 &&\hspace{2cm}-
  E^{\ell,p}\left(\prod_{1\leq i\leq
p}(\vr_i-s)^{q-p}I_{\vr_i\geq s}\right)
  \\
  &\leq&
  \sum_{{\tiny\begin{array}{l}
       0\leq k_i\leq q-p \\
       1\leq i\leq p \\
       \sum k_i\geq 1
       \end{array}}}
    \left(\frac{q}{\sqrt{2x}}\right)^{\sum k_i}
    \prod_{i=1}^{p}\left({q-p}\atop {k_i}\right)
   E^{\ell,p} \prod_{i=1}^p(\vr_i-s)^{q-p-k_i}
   I_{\vr_i\geq s}\eean
Note that this estimate holds for all $\ell$, with all
expressions vanishing when $\ell < p(n-p)$. Using the
previous estimate and using (5.0.1), one finds the
estimate
 \bean 0&\leq&{
\frac{1}{(2x)^{\frac{p(q-p)}{2}}}E_{x,p}\left(
I_{\lb\supseteq\mu}(\lb)\prod_{(i,j)\in\lb\atop{n-q<j\leq
n-p}} h^{\lb}_{(i,j)}\right)
 }\\
&& \hspace{2cm}-
 e^{-px}\sum_{\ell\geq 0}
\frac{(px)^{\ell}}{\ell !} E^{\ell,p}\left( \prod_{1\leq
i\leq p}(\vr_i-s)^{q-p}I_{\vr_i\geq s}
\right) % \\
\eean\bean &\leq&
  \sum_{{\footnotesize\begin{array}{l}
       0\leq k_i\leq q-p \\
       1\leq i\leq p \\
       \sum k_i\geq 1
       \end{array}}}
    \left(\frac{q}{\sqrt{2x}}\right)^{\sum k_i}
    \prod_{i=1}^{p}\left({q-p}\atop {k_i}\right)
\\&&\hspace{2cm}
    E_{x,p}\left(
\prod_{i=1}^p(\frac{\lb_i-x}{\sqrt{2x}}-s)^{q-p-k_i}
%I_{\vr_i\geq c}
I_{\frac{\lb_i-x}{\sqrt{2x}}\geq s}\right)
% \eean\bean%
 \\
 &=& O(\frac{1}{\sqrt{2x}})
 ~~~\mbox{for $x\rightarrow \iy$}
   \eean
because the expectation $E_{x,p}$ tends to an integral,
by Johansson's theorem, applied to the function
 $$
g(x_1,\ldots,x_p)=\prod_1^p(x_i-s)^{q-p-k_i}I_{[s,\iy)}(x_i).
$$
%
%\newpage
%
For $0\leq\ell <p(n-p)$, the expectation is
automatically zero.
Then
\bean
\lefteqn{\lim_{x\rg\iy}\frac{1}{(2x)^{\frac{p(q-p)}{2}}}E_{x,p}\left(
I_{\lb\supseteq\mu}(\lb)\prod_{(i,j)\in\lb\atop{n-q<j\leq n-p}}
h^{\lb}_{(i,j)}\right)}\\
 & & \\
%&=&\lim_{x\rg\iy}\frac{e^{-px}}{(2x)^{\frac{p(q-p)}{2}}}\sum_{\ell\geq
%p(n-p)}\frac{(px)^{\ell}}{\ell !}E^{\ell,p}\left(
%I_{\lb\supseteq\mu}(\lb)\prod_{(i,j)\in\lb\atop{n-q<j\leq n-p}}
%h^{\lb}_{(i,j)}\right)\\
% & & \\
%& &\hspace*{5cm}\mbox{(notice that $\ell\geq p(n-p)$ in the}\\
%& &\hspace*{5cm}\mbox{summation is no restriction)}\\
 & & \\
 &=&\lim_{x\rg\iy}e^{-px}\sum_{\ell\geq 0}
\frac{(px)^{\ell}}{\ell !} E^{\ell,p}\left( \prod_{1\leq
i\leq p}(\vr_i-s)^{q-p}I_{\vr_i\geq s}
\right)  \\
 & & \\
 &=&\lim_{x\rg\iy}e^{-px}\sum_{\ell\geq 0}
\frac{(px)^{\ell}}{\ell !}E^{\ell,p}\left( \prod_{1\leq
i\leq p}\left(\frac{\lb_i-x}{\sqrt{2x}}-s\right)^{q-p}
I_{\frac{\lb_i-x}{\sqrt{2x}}\geq s}\right) \\
& &\\
&=&\lim_{x\rg\iy}E_{x,p}\left( \prod_{1\leq i\leq
p}\left(\frac{\lb_i-x}{\sqrt{2x}}-s\right)^{q-p}
I_{\frac{\lb_i-x}{\sqrt{2x}}\geq s}\right) \\
& & \\
&=&p!\int_{x_1>\ldots >x_p}\prod_{1\leq i\leq
p}\left((x_i-s)^{q-p}I_{[s,\iy)}(x_i)\right)
 \frac{1}{Z_p}\Delta_p(x)^2\prod^p_{j=1}e^{-x^2_j}
%\vp_p(x_1,\ldots,x_p)
dx_j  ,
 \eean
   by applying Proposition 5.1 to the continuous
function
$$
g(x_1,\ldots,x_p)=\prod_1^p(x_i-s)^{q-p}I_{[s,\iy)}(x_i).
$$
This leads to the ratio of Hermitian matrix integrals in
(\ref{main theorem}). The numerator is precisely the
matrix integral (\ref{matrix integral}), which according
to Theorem 3.1 satisfies the Painlev\'e equation
(\ref{ode in h}), which ends the proof of Theorem
1.4.\qed

%\newpage

\section{Painlev\'e IV as a limit of Painlev\'e V}

According to Theorem 1.4, the limit (\ref{main theorem})
lead to a solution of Painlev\'e IV; this was
established by identifying the limit as a matrix
integral and then applying Theorem 3.1. In this section
we give an alternative proof of this fact, by directly
taking the scaling limit of the Painlev\'e V equation
(\ref{1.0.14}).
Remember from Proposition 1.1 and Theorem 1.4, we have
\bea
 \lefteqn{\prod_{i=1}^p \frac{(n-i)!}{(q-i)!}
  x^{-(n-p)p}e^{px }
   ~E_{x,p}\left( I_{\{\lambda \supseteq \mu \}}(\lb)
  \prod_{{(i,j)\in \lb}
 \atop {n-q<j\leq n-p}}
  h_{(i,j)}^{\lb}
   \right)}\no\\
 &&~~~~~~~~~~~~~~~~~~~~~~~~=
%    \tilde{\tilde c_{ }}
%    x^{(n-p)p}
    \exp {\displaystyle{\int_0^x\frac{u(y)-p(n-p)+py}{y}
    dy}
    }
\label{6.0.1} \eea and, with the rescaling
\be
 n-p=x+\sqrt{2x}s,
 \label{rescaling1}\ee
we have from section 5, that
 \be
 \lim_{x\rg\iy}\frac{1}{(2x)^{\frac{p(q-p)}{2}}}
E_{x,p}\left(I_{\lb\supseteq\mu}(\lb)\prod_{(i,j)\in
\lb\atop{n-q<j\leq n-p}} h^{\lb}_{ij}\right)=
 c ~\exp{\int_0^s h(y)dy} %\label{main theorem}
 .\label{6.0.3}\ee
 for an appropriate choice of $h(y)$ and $c$. Taking into account the inverse of the
 rescaling (\ref{rescaling1}),
 which for large $n$ reads
  \be
   x=n-s\sqrt{2n} +o(1)\label{approximate rescaling}
   \ee
   and taking the logarithmic derivatives of both equations
   (\ref{6.0.1}) and (\ref{6.0.3}), we have
   \bean
     \frac{\pl}{\pl s} \log E_{x,p}
      \left( I_{\{\lambda \supseteq \mu \}}(\lb)
  \prod_{{(i,j)\in \lb}
 \atop {n-q<j\leq n-p}}
  h_{(i,j)}^{\lb}
   \right)
  &=&
    \frac{u(x)}{x}\frac{\pl x}{\pl s}
%   \\ &simeq$
   \\
   \frac{p(q-p)}{2x} \frac{\pl x}{\pl s} +
   \frac{\pl}{\pl s} \log E_{x,p}
      \left( I_{\{\lambda \supseteq \mu \}}(\lb)
  \prod_{{(i,j)\in \lb}
 \atop {n-q<j\leq n-p}}
  h_{(i,j)}^{\lb}
   \right) & \simeq&  h(s),
  \\ &&\hspace{-2cm}\mbox{for large, but fixed $n$.}
   \eean
Letting $n\rightarrow \iy$ in both equations and keeping
the leading terms lead to

\bean
 \frac{\pl}{\pl s} \log E_{x,p}
    & \simeq&   -\sqrt{\frac{2}{n}}u(n-s\sqrt{2n} )\\
  \frac{\pl}{\pl s} \log E_{x,p}
    & \simeq& h(s)
    .\eean
 So, setting
  $$
  h(s):= -\sqrt{\frac{2}{n}}u(n-s\sqrt{2n} )
  ,$$
  we have
  $$
  u(x)
  =-\sqrt{\frac{n}{2}}
   h\left(\frac{n-x}{\sqrt{2n}}\right)
  $$
  and thus
  $$
%  \begin{array}{ll}
   u'(x)=\frac{1}{2} h'~,\qquad
   u''(x)=-\frac{1}{2\sqrt{2n}} h''
   ~,\qquad
   u'''(x)=-\frac{1}{4n} h'''
.$$
Setting the rescaling (\ref{approximate rescaling}) in
the polynomials $Q(x)$ and $R(x)$, as defined in
(\ref{QR-polynomials}), we have
 \bean
  4Q&=& n\left(-2 s^2+4(2p-q) \right)+o(n)\\
  4Q'&=& 2\sqrt{2n}s+o(\sqrt{n})  \\
  2R&=& 2np(p-q)+o(n).%\\
  \eean
Substituting into (\ref{1.0.14}), keeping the leading
terms, which are of order $n$, and multiplying by $4$,
one finds the equation (\ref{ode in h-1}),
\be
 h^{\prime\prime\prime}+6h^{\prime
2}-4(y^2+2(q-2p))h^{\prime}+4yh-8(q-p)p =0.
%\label{ode in h-1}
 \ee

\section{Appendix: Chazy classes}

In his classification of differential equations
$$
f^{\prime\prime\prime}=F(z,f,f^{\prime},f^{\prime\prime}),~\mbox{where
$F$ is rational in $f,f^{\prime},f^{\prime\prime}$ and
locally analytic in z,}$$ subjected to the requirement
that the general solution be free of movable branch
points,
 Chazy found thirteen cases, the first
of which is given by \be f^{\prime
\prime\prime}+\frac{P'}{P}f^{
\prime\prime}+\frac{6}{P}f^{\prime 2}-\frac{4P'}{P^2}ff'
+\frac{P^{\prime\prime}}{P^2} f^2
+\frac{4Q}{P^2}f'-\frac{2Q'}{P^2}f+\frac{2R}{P^2}=0 \ee
 with arbitrary polynomials $P(z), Q(z), R(z)$ of
 maximal
degree $3,2,1$ respectively. Cosgrove and Scoufis
\cite{CosgroveScoufis,Cosgrove}, (A.3), show that this
third order equation
 has a first integral, which is second order in $f$
and quadratic in $f^{\prime\prime}$, \bea &f^{
\prime\prime 2}& +\frac{4}{P^2}
 \left( (Pf^{\prime 2}+Q f^{\prime}+R)f^{\prime}
 - (P' f^{\prime 2}+\frac{}{}Q' f^{\prime}+R')f^{}
  \right. \nonumber\\ && \hspace{2cm} \left.
   +\frac{1}{2}(P^{\prime\prime}f^{\prime
}+Q^{\prime\prime} )f^2
 -\frac{1}{6} P^{\prime\prime\prime}f^3 +\dt\right)=0;
\label{Cosgrove}\eea $\dt$ is the integration constant.

 $$
 \begin{array}{c|c|c|c} P&4Q&2R&\mbox{Painlev\'e eqt}\\ \hline\\
 1 &  -4(u^2 + 2(a-n))&  -8an & \mbox{P}\,IV
 \\
 u &  -u^2+2u(2n+a-b)-(a+b)^2 &  -bn(a+b+u)& \mbox{P}\,V.
 \end{array}$$

%\newpage

 Equations of the
general form $$ f^{ \prime\prime 2}=G(x,f,f^{ \prime})$$
are invariant under the map $$ x\mapsto
\frac{a_1z+a_2}{a_3z+a_4}~~ \mbox{and}~~ f\mapsto
\frac{a_5f+a_6z+a_7}{a_3z+a_4}.$$ Using this map, the
polynomial $P(z)$ can be normalized to $$
P(z)=z(z-1),~z,~\mbox{or} ~1.$$

 In this way, Cosgrove shows (7.0.2) is a master Painlev\'e
 equation, containing the 6 Painlev\'e equations.  In
the cases of PIV and PV, the canonical equations are
respectively:

$$  g^{\prime\prime 2}=-4g^{\prime 3
}+4(zg^{\prime}-g)^2+A_1g^{\prime}+A_2
\hspace{1cm}(\mbox{\bf Painlev\'e IV})$$

$$ (zg^{\prime\prime
})^2=(zg^{\prime}-g)\Bigl(-4g^{\prime 2
}+A_1(zg^{\prime}-g)+A_2\Bigr)+A_3g^{\prime}+A_4
.%
 \hspace{1cm}(\mbox{\bf
Painlev\'e V})
 $$

\end{document}